\documentclass[a4paper,12pt]{elsarticle}
\usepackage{amsmath,amssymb,theorem}
\usepackage{graphicx}
\usepackage{subfigure}
\usepackage{psfrag}
\usepackage{color}
\usepackage{epstopdf}
\usepackage[T1]{fontenc}

\newtheorem{thm}{Theorem}[section]

\newtheorem{obs}[thm]{Observation}
\newtheorem{cor}[thm]{Corollary}
\newtheorem{lem}[thm]{Lemma}
\theorembodyfont{\rmfamily}
\newtheorem{rem}[thm]{Remark}
\def\pf{\bigskip\noindent {\bf Proof.}~~}

\def\dfn#1{{\sl #1}}

\def\less{\backslash}

\def\pf{\bigskip\noindent {\emph{Proof.}}~~}

\def\qed{ \hfill $\square$}

\newcounter{counter}

\def\proofsquare{  \bigskip\hfill\vrule height3pt width6pt depth2pt}

\begin{document}
\begin{frontmatter}
\title{$k$-Ary spanning trees contained in tournaments}

\author[1,2]{Jiangdong Ai}
\ead{Jiangdong.Ai.2018@live.rhul.ac.uk}

\author[3]{Hui Lei\corref{cor1}}
\ead{hlei@nankai.edu.cn}

\author[2]{Yongtang Shi}
\ead{shi@nankai.edu.cn}

\author[2]{Shunyu Yao}
\ead{phpass@mail.nankai.edu.cn}

\author[4]{Zan-Bo Zhang}
\ead{eltonzhang2001@gmail.com}

\address[1]{Department of Computer Science,
Royal Holloway, University of London,
Egham, Surrey, TW20 0EX, UK}
\address[2]{Center for Combinatorics and LPMC, Nankai University, Tianjin 300071, China}
\address[3]{School of Statistics and Data Science, LPMC and KLMDASR, Nankai University, Tianjin 300071, China}
\address[4]{Institute of Artificial Intelligence \& Deep Learning, and School of Statistics \& Mathematics, Guangdong University of Finance \& Economics, Guangzhou, 510320, China}
\cortext[cor1]{Corresponding author}

\begin{abstract}

A rooted tree is called a \dfn{$k$-ary tree}, if all non-leaf vertices have exactly $k$ children, except possibly one non-leaf vertex has at most $k-1$ children. Denote by $h(k)$ the minimum integer such that every tournament of order at least $h(k)$ contains a $k$-ary spanning tree. It is well-known that every tournament contains a Hamiltonian path, which implies that $h(1)=1$. Lu et al. [J. Graph Theory {\bf 30}(1999) 167--176] proved the existence of $h(k)$, and showed that $h(2)=4$ and $h(3)=8$. The exact values of $h(k)$ remain unknown for $k\geq 4$. A result of Erd\H{o}s on the domination number of tournaments implies $h(k)=\Omega(k\log k)$. In this paper, we prove that $h(4)=10$ and $h(5)\geq13$.\\

\end{abstract}
\end{frontmatter}

\section{Introduction}

\baselineskip 18pt
In this paper, we consider digraphs which are finite and simple. That is, we do not permit the existence of loops or multiple directed arcs. For any undefined terms about digraphs, we refer the reader to the book of Bang-Jensen and Gutin \cite{BG}.

A {\it tournament} $T=(V,E)$ is a directed graph (digraph) obtained by assigning a direction for each edge in an undirected complete graph. In this paper, a tournament of order $n$ is called \dfn{$n$-tournament}. We also use $x\to y$ or $(x,y)$ to denote an arc $xy\in E$, say $x$ {\it beats} $y$. Let $A\Rightarrow B$ denote that every vertex in $A$ beats every vertex in $B$. We call a tournament {\it transitive} if $x\to y$ and $y\to z$ imply that $x\to z$, in other words, its vertices can be linearly ordered such that each vertex beats all later vertices.  We denote by $[X_i]$ the vertex set $\{x_1,\ldots,x_i\}$ for $i\geq1$.
If $x\to y$, we call $y$ an \dfn{out-neighbor} of $x$, and $x$ an \dfn{in-neighbor} of $y$. We use $N^+(x)$ and $N^-(x)$ to denote the \dfn{out-neighborhood} and the \dfn{in-neighborhood} of a vertex $x$ of $T$, respectively. Correspondingly, we use $d^+(x)=|N^+(x)|$ and $d^-(x)=|N^-(x)|$ to denote the \dfn{out-degree} and the \dfn{in-degree} of a vertex $x$ of $T$, respectively. A {\it leaf} is a vertex of  out-degree zero.
For $x\in V$ and $X\subseteq V$, we denote by $N^+_{X}(x)$ (resp. $N^-_{X}(x)$) the set of out(resp. in)-neighborhood of $x$ in $X$, that is, $N^+_{X}(x)=N^+(x)\cap X$ (resp. $N^-_{X}(x)=N^-(x)\cap X$) (here, $x$ may or may not belong to $X$).
We write $d^+_X(x) = \left|N^+_{X}(x)\right| $, $d^-_X(x) = \left|N^-_{X}(x)\right| $ and $d^+_X=\max\{d^+_X(v)|v\in X\}$.
A tournament is \dfn{$k$-regular} if all vertices have in-degree and  out-degree  $k$.
For a subset $X\subseteq V$, we denote by $T[X]$ the subtournament of $T$ induced by $X$.

A {\it rooted tree} is a directed tree with a special vertex, called the {\it root}, such that there exists a unique (directed) path from the root to any other vertex. A rooted tree is called a {\it $k$-ary tree}, if all non-leaf vertices have exactly $k$ children, except possibly one non-leaf vertex has at most $k-1$ children.  If all non-leaf vertices have exactly $k$ children, then we call it a {\it full $k$-ary tree}. A \dfn{$k$-star} is a full $(k-1)$-ary tree with $k$ vertices.

An oriented graph $H$ on $n$ vertices is {\it unavoidable} if every $n$-tournament contains $H$ as a subgraph, otherwise, we say that $H$ is {\it avoidable}. The concept of unavoidable was introduced by Linial et al. \cite{LSS}, in which they studied the maximum number of edges that an unavoidable subgraph on $n$ vertices can have.
In particular, if $H$ contains a directed cycle then $H$ must be avoidable, since a transitive tournament contains no directed cycles and hence no copy of $H$. It is therefore natural to ask which oriented trees are unavoidable.

R\'edei \cite{R} showed that every tournament contains a Hamiltonian path. Thomason \cite{5} proved that all orientations of sufficiently long cycles are unavoidable except for those which yield directed cycles. Erd\H{o}s \cite{R1989}  proved that for any fixed positive integer $m$, there exists a number $f(m)$ such that every $n$-tournament contains $\left \lfloor \frac{n}{m}\right \rfloor$ vertex-disjoint transitive sub-tournaments of order $m$ if $n\ge f(m)$. H\"aggkvist and Thomason \cite{HT} showed that every oriented tree of order $m$ is contained in every tournament of order $12m$ and El Sahili \cite{El} improved the bound to $3(m-1)$. Lu et al. \cite{L1993, LWW1998} investigated the avoidable claws.
For more results on unavoidable digraphs, we refer to \cite{G1971,L1996,R1974}.

Actually, R\'edei's result \cite{R} can be restated as that a $1$-ary spanning tree is unavoidable. It is therefore natural to study the general problem of whether a tournament contains a $k$-ary spanning tree.
Lu et al. \cite{LWCLW} proved the following fundamental theorem for the existence of a $k$-ary spanning tree of a tournament.
\begin{thm}[\cite{LWCLW}]\label{existence}
For any fixed positive integer $k$, there exists a number $h'(k)$ such that every $n$-tournament contains a $k$-ary spanning tree if $n\geq h'(k)$.
\end{thm}

Define $h(k)$ as the minimum number such that every tournament of order at least $h(k)$ contains a $k$-ary spanning tree. The existence of a Hamiltonian
path for any tournament is the same as $h(1) = 1$.
Lu et al. \cite{LWCLW} determined that $h(2)=4$ and $h(3)=8$. The exact values of $h(k)$ remain unknown for $k\geq 4$.
A result of Erd\H{o}s on the domination number of tournaments implies $h(k)=\Omega(k\log k)$.
\begin{thm}\label{lowerbound}
For any $k\geq4$, $h(k)=\Omega(k\log k)$.
\end{thm}

In this paper, we prove that $h(4)=10$ and $h(5)\geq13$.
\begin{thm}\label{main}
$h(4)=10$ and $h(5)\ge 13$.
\end{thm}



\section{Proof of Theorem \ref{lowerbound}}
For any $X,Y\subseteq V(T)$, we say that $X$ {\it dominates} $Y$ if for every $v\in Y\less X$ there exists a $u\in X$ which beats $v$. The {\it domination number} of $T$, denoted $\mu(T)$, is the smallest cardinality of a set that dominates $V(T)$.

Erd\H{o}s \cite{E1963} used the probabilistic method to prove the following fact.
\begin{lem}[\cite{E1963}]\label{pm}
For every $\varepsilon>0$ there is a number $K$ such that for for every $k\geq K$ there exists a tournament $T_k$ with no more than $2^kk^2\log(2+\varepsilon)$ vertices such that $\mu(T_k)>k$.
\end{lem}
By Lemma \ref{pm}, we can get the following Corollary \ref{donu} directly which is stated in \cite{MV}.
\begin{cor}[\cite{MV}]\label{donu}
There exists a constant $c>0$ such that for every $n$ there exists a tournament $T$ with $n$ vertices such that $\mu(T)>c\log n$.
\end{cor}

Now we present the proof of Theorem \ref{lowerbound}.\\

\noindent {\it Proof of Theorem \ref{lowerbound}.} \\
Let $T$ be a tournament with $n$ vertices and $\mu(T)>c\log n$. Suppose $T$ contains a $k$-ary spanning tree $R$. Since the number of the non-leaf vertices of $R$ is $\lceil\frac{n-1}{k}\rceil$ and all non-leaf vertices of $R$ dominates $V(T)$, we have $\lceil\frac{n-1}{k}\rceil\geq\mu(T)$. Then $n> (\mu(T)-1)k+1$. By Corollary \ref{donu}, we have $h(k)=\Omega(k\log k)$.\qed

\section{Proof of Theorem \ref{main}}
We need the following three useful lemmas proved in \cite{LWCLW}.
\begin{lem}[\cite{LWCLW}]\label{induction}
Let $R$ be a $k$-ary tree of tournament $T$ with the root $v$ and $S$ a $k$-star of $T$ with the root $u$, where $R$ and $S$ are vertex disjoint. If $d_{V(R)}^+(u)\ge 1$, then $T$ contains a $k$-ary tree $R'$ with $V (R') = V (R) \cup V (S)$.
Furthermore, if $u \in N^+(v)$, then $R'$ can be chosen to have the root $v$, which is the same root as $R$.
\end{lem}

\begin{lem}[\cite{LWCLW}]\label{km}
If every $(km + 1)$-tournament has a $k$-ary spanning tree, then so does every $km$-tournament.
\end{lem}

According to the structure of  $k$-ary spanning trees, we can directly obtain the following result.
\begin{obs}\label{obs}
For any $n$-tournament $T=(V,E)$ with $n\geq 2k+1$, let $T_{\geq k}=\{v\in V~|~d^+(v)\geq k\}$. If for any two different vertices $u,v\in T_{\geq k}$, $|(N^+(u)\cup N^+(v))\less \{u,v\}|\leq 2k-2$, then $T$ contains no $k$-ary spanning tree.
\end{obs}



\noindent {\it Proof of Theorem \ref{main}.} \\
First, we consider the case of $k=4$. Let $T_9$ be the $9$-tournament with $V(T_9)=\{0,1,\ldots,8\}$ and $E(T_9)=\{ij:i-j\equiv1,2,3,5 \ (\text{mod} \ 9)\}$. By Observation \ref{obs}, it is straightforward to check that $T_9$ does not contain a $4$-ary spanning tree, since $d_{V(T_9)}^+(i)=4$ for any $i\in V(T_9)$, and $N^+_{V(T_9)}(j)\cap N^+_{V(T_9)}(i)\neq \emptyset$ for any $j\in N^+_{V(T_9)}(i)$. So $h(4)\geq10$. In the following, by induction, we will prove that every tournament $T$ of order $n\geq 10$ contains a $4$-ary spanning tree.

Let $T=(V,E)$ be a tournament of order $n$. Note that for any $X\subseteq V$, we have $d_X^+\geq\left\lceil\frac{|X|-1}{2}\right\rceil$. Suppose $n\geq14$ and the theorem is true for all $n'<n$. Since
$n\geq14$, we can choose $v\in V$ with $d^+(v)\geq4$, say $N^+(v)=\{a,b,c,d\}$. Let $T'=T[V\less\{v,a,b,c\}]$. By the induction hypothesis, $T'$ contains a $4$-ary spanning tree. By Lemma \ref{induction}, $T$ contains a $4$-ary spanning tree. Therefore, by Lemma \ref{km}, it suffices to prove that every tournament $T$ of order $n$ contains a $4$-ary spanning tree, where $n\in\{10,11,13\}$. Let $u$ be a vertex of $T$ with the maximum out-degree and $V=\{u\}\cup[X_{n-1}]$.\\

\noindent {\bf Claim\refstepcounter{counter}\label{234}  \arabic{counter}.} For $1\leq d^-(u)\leq4$, if there exists a vertex $v\in N^-(u)$ such that $d_{N^-(u)}^+(v)=d^-(u)-1$ and $d_{N^+(u)}^+(v)\geq 4-d^-(u)$,
then $T$ contains a $4$-ary spanning tree.

\pf Let $N^-(u)=[X_{d^-(u)}]$. Suppose $d_{N^-(u)}^+(x_1)=d^-(u)-1$ and $d_{N^+(u)}^+(x_1)\geq 4-d^-(u)$, say $x_1\Rightarrow\{x_2,x_3,x_4\}$.
Since $d_{\{x_5,\dots,x_{n-1}\}}^+\ge \lceil\frac{n-6}{2}\rceil\geq n-9$, we may assume $x_8\Rightarrow\{x_9,\dots,x_{n-1}\}$. Then we obtain a $4$-ary spanning tree of $T$ induced by $\{x_1x_2,x_1x_3,x_1x_4,x_1u,ux_5,\dots,ux_8,\\ x_8x_9,\dots,x_8x_{n-1}\}$.
\proofsquare

\noindent {\bf Claim\refstepcounter{counter}\label{10ve}  \arabic{counter}.}
Every $10$-tournament $T$ contains a $4$-ary spanning tree.

\pf We consider the following five cases.\medskip

{\bf Case 1:} $d^+(u)=9$.

Since $d^+_{[X_9]}\geq4$, we assume $x_9\Rightarrow[X_4]$ and $x_6\to x_5$. Then we obtain a $4$-ary spanning tree induced by $\{ux_6,\dots,ux_9,x_9x_1,\dots,x_9x_4,x_6x_5\}$.\medskip

{\bf Case 2:} $d^+(u)=8$, say $N^+(u)=[X_8]$.

Since $d^+_{[X_8]}\geq4$, we assume $x_8\Rightarrow[X_4]$. Then we obtain a $4$-ary spanning tree induced by $\{ux_5,\dots,ux_8,x_8x_1,\dots,x_8x_4,x_9u\}$.\medskip

{\bf Case 3:} $d^+(u)=7$, say $N^+(u)=[X_7]$ and $x_9\to x_8$.

By Claim \ref{234}, we may assume
$d_{[X_7]}^+(x_9)\le1$. If $d_{[X_7]}^+(x_8)\geq3$,
assume $x_8\Rightarrow\{x_7,x_6,x_5\}$,
and then $\{x_9x_8,x_8x_5,x_8x_6,x_8x_7,x_8u,ux_1,\dots,ux_4\}$ induces a desired $4$-ary spanning tree.
So we may assume $d_{[X_7]}^+(x_8)\le2$. Then $\big|N_{[X_7]}^-(x_9)\cap N_{[X_7]}^-(x_8)\big|\geq4$,
say $[X_4]\Rightarrow\{x_8,x_9\}$. Without loss of generality, we may assume that $d_{[X_4]}^+(x_3)\ge2$ with $x_3\Rightarrow[X_2]$ and $x_6\to x_7$. Then we obtain a $4$-ary spanning tree induced by $\{ux_3,\dots,ux_6, x_3x_1,x_3x_2, x_3x_8,x_3x_9,x_6x_7\}$.\medskip

{\bf Case 4:} $d^+(u)=6$, say $N^+(u)=[X_6]$ and $x_9\rightarrow x_8,x_8\rightarrow x_7$.

If $d_{[X_6]}^+(x_9)\ge2$ or $d_{[X_6]}^+(x_8)\ge2$,
say $x_9\Rightarrow\{x_5,x_6\}$ or $x_8\Rightarrow\{x_5,x_6\}$, then we obtain a $4$-ary spanning tree induced by $\{x_9x_5,x_9x_6,x_9x_8,x_9u,ux_1,\dots,ux_4,x_8x_7\}$ or $\{x_9x_8,x_8x_5,x_8x_6,\\ x_8x_7,x_8u,ux_1,\dots,ux_4\}$. So we assume $\big|N_{[X_6]}^-(x_9)\cap N_{[X_6]}^-(x_8)\big|\geq4$, say $[X_4]\Rightarrow\{x_9,x_8\}$. Without loss of generality, we may assume that $d_{[X_4]}^+(x_3)\ge2$ with $x_3\Rightarrow[X_2]$. Then we obtain a desired $4$-ary spanning tree induced by $\{ux_3,\dots,ux_6,x_3x_1,x_3x_2,x_3x_8,x_3x_9,x_8x_7\}$.\medskip

{\bf Case 5:} $d^+(u)=5$, say $N^+(u)=[X_5]$.

By Claim \ref{234}, we may assume $d_{N^-(u)}^+\leq2$. Without loss of generality, we may assume that $x_9\Rightarrow\{x_7,x_8\}$, $x_8\Rightarrow\{x_7,x_6\}$ and $x_6\to x_9$. If $d_{[X_5]}^+(x_9)\ge1$ or $d_{[X_5]}^+(x_8)\ge1$,
say $x_9\to x_5$ or $x_8\to x_5$, then one can find a desired tree induced by $\{x_6x_9,x_9x_5,x_9x_7,x_9x_8,x_9u,ux_1,\dots,ux_4\}$ or $\{x_9x_8,x_8x_5,x_8x_6,x_8x_7,x_8u,ux_1,\dots,ux_4\}$. So we may further assume $[X_5]\Rightarrow\{x_9,x_8\}$. Without loss of generality, assume that $x_7\to x_6$. Since $d^+(x_7)\leq5$, we have $\big|N_{[X_5]}^-(x_7)\big|\geq2$, say $\{x_4,x_5\}\Rightarrow x_7$ and $x_4\to x_5$.
Then the set $\{ux_1,\dots,ux_4,x_4x_5,x_4x_9,x_4x_8,x_4x_7,x_7x_6\}$ induces a desired $4$-ary spanning tree.
\proofsquare

Suppose $n\in\{11,13\}$ and
$d^+(u)=n-1$. By Claim \ref{10ve}, let $R$ be a $4$-ary spanning tree of $T[[X_{10}]]$.  Without loss of generality, we assume that $R'\subseteq R$ is a full $4$-ary tree rooted at $x_9$ with $V(R')=[X_{9}]$. Then $R'\cup\{ux_9,\dots,ux_{n-1}\}$ induces a desired $4$-ary spanning tree. So we may further assume $n\in\{11,13\}$
and $N^+(u)=[X_{d^+(u)}]$ with $d^+(u)\leq n-2$ in the following.\\

\noindent {\bf Claim\refstepcounter{counter}\label{11ve}  \arabic{counter}.}
Every $11$-tournament $T$ contains a $4$-ary spanning tree.

\pf We consider the following five cases.\medskip

{\bf Case 1:} $d^+(u)=9$.

By Claim \ref{234}, we may assume
$[X_7]\Rightarrow x_{10}$. Without loss of generality, we assume $x_4\Rightarrow[X_3]$ since $d_{[X_7]}^+\ge3$, and $x_7\Rightarrow\{x_8,x_9\}$ since $d_{\{x_5,\dots,x_9\}}^+\ge2$. Then we obtain a $4$-ary spanning tree of $T$ induced by $\{ux_4,\dots,ux_7,x_4x_1,x_4x_2,x_4x_3,x_4x_{10},x_7x_8,x_7x_9\}$.\medskip

{\bf Case 2:} $d^+(u)=8$.

Without loss of generality, we assume that $x_{10}\rightarrow x_9$ and $x_4\Rightarrow\{x_5,x_6,x_7,x_8\}$ because $d_{[X_8]}^+\ge4$. Then we find a desired $4$-ary spanning tree induced by $\{x_{10}x_9,x_{10}u, ux_1,\dots,ux_4,x_4x_5,\\ \dots,x_4x_8\}$.\medskip

{\bf Case 3:} $d^+(u)=7$.

Suppose $x_{10}\Rightarrow \{x_8,x_9\}$. By Claim \ref{234}, we may assume $[X_7]\Rightarrow x_{10}$ and $x_4\Rightarrow \{x_5,x_6,x_7\}$ since $d_{[X_7]}^+\ge 3$. We obtain a $4$-ary spanning tree of $T$ induced by $\{ux_1,\dots,ux_4,x_4x_5,x_4x_6,x_4x_7,\\x_4x_{10},x_{10}x_9,x_{10}x_8\}$. Suppose $x_{10}\to x_9$, $x_9\to x_8$ and $x_8\to x_{10}$. If $d_{[X_7]}^+(x_9)\geq3$, say $x_9\Rightarrow\{x_5,x_6,x_7\}$, then we obtain a $4$-ary spanning tree of $T$ induced by $\{x_{10}x_9,x_{10}u, ux_1,\dots,ux_4,\\x_9x_5,\dots,x_9x_8\}$. If $x_9\Rightarrow\{x_6,x_7\}$ and $x_8\to x_5$, then we obtain a desired $4$-ary spanning tree induced by $\{x_9x_6,x_9x_7,x_9x_8,x_9u, ux_1,\dots,ux_4,x_8x_{10},x_8x_5\}$. By the symmetry of $x_8$ and $x_9$, we may assume $[X_4]\Rightarrow\{x_8,x_9\}$ and $x_1\Rightarrow\{x_2,x_3\}$ since $d^+_{[X_4]}\geq2$. Then $\{x_1x_2,x_1x_3,x_1x_8,x_1x_9,x_8x_{10},x_8u,ux_4,\dots,ux_7\}$ induces a desired $4$-ary spanning tree.\medskip

{\bf Case 4:}  $d^+(u)=6$.

By Claim \ref{234}, we may assume $d_{N^-(u)}^+\leq2$. Let $x_{10}\Rightarrow\{x_9,x_8\}$, $x_9\Rightarrow\{x_8,x_7\}$ and $x_7\to x_{10}$. If $d_{[X_6]}^+(x_{10})\geq2$ or $d_{[X_6]}^+(x_9)\geq2$, say $x_{10}\Rightarrow\{x_5,x_6\}$ or $x_9\Rightarrow\{x_5,x_6\}$, then we obtain a desired tree induced by $\{x_7x_{10},x_7u,x_{10}x_9,x_{10}x_8,x_{10}x_6,x_{10}x_5,ux_1,\dots,ux_4\}$ or $\{x_{10}x_9,x_{10}u,x_9x_5,\dots,x_9x_8,ux_1,\dots,ux_4\}$.
So we may further assume $[X_4]\Rightarrow\{x_{10},x_9\}$ and $x_3\Rightarrow[X_2]$ because $d_{[X_4]}^+\ge2$. Then we obtain a $4$-ary spanning tree induced by $\{ux_3,\dots,ux_6,x_3x_1,x_3x_2,x_3x_{10},x_3x_9,x_9x_8,x_9x_7\}$.\medskip

{\bf Case 5:} $d^+(u)=5$.

In this case, $T$ is a $5$-regular tournament. Let $1\le d^+_{[X_5]}(x_4)\le 2$ with $x_4\to x_5$. We may assume  $x_4\Rightarrow \{x_6,x_7,x_8\}$ because $d^+(x_4)=5$, and let $x_{10}\to x_9$. Then we obtain a $4$-ary spanning tree induced by $\{x_{10}x_9,x_{10}u,ux_1,\dots,ux_4,x_4x_5,\dots,x_4x_8\}$.
\proofsquare

\noindent {\bf Claim\refstepcounter{counter}\label{13ve}  \arabic{counter}.}
Every $13$-tournament $T$ contains a $4$-ary spanning tree.

\pf We consider the following six cases.\medskip

{\bf Case 1:} $d^+(u)=11$.

By Claim \ref{234}, we may assume $[X_9]\Rightarrow x_{12}$ and $x_4\Rightarrow[X_3]$ because $d_{[X_9]}^+\geq4$. First we suppose $d_{\{x_5,\dots,x_{11}\}}^+\geq4$, say
$x_7\Rightarrow \{x_8,\dots,x_{11}\}$. Then we obtain a $4$-ary spanning tree induced by $\{ux_4,\dots,ux_7,x_7x_8,\dots,x_7x_{11},x_4x_1,x_4x_2,x_4x_3,x_4x_{12}\}$. Next we consider $d_{\{x_5,\dots,x_{11}\}}^+=3$.
If $x_4\Rightarrow\{x_5,\dots,x_{11}\}$, then $d^+(x_4)=11$. Since $d_{N^+(x_4)}^+(u)\geq 3$, we obtain a $4$-ary spanning tree of $T$ by Claim \ref{234}. Otherwise there exists a vertex $v\in \{x_5,\dots,x_{11}\}$ such that $v\to x_4$, without loss of generality, say $x_8\Rightarrow \{x_4,\dots,x_7\}$. Then we obtain a $4$-ary spanning tree induced by $\{ux_8,\dots,ux_{11},x_8x_4,\dots,x_8x_7,x_4x_1,x_4x_2,x_4x_3,x_4x_{12}\}$.\medskip

{\bf Case 2:} $d^+(u)=10$, say $x_{12}\to x_{11}$.

By Claim \ref{234}, we may assume $[X_9]\Rightarrow x_{12}$. If $d_{[X_9]}^+(x_{11})\geq3$, say $x_{11}\Rightarrow\{x_7,x_8,x_9\}$ and $x_3\Rightarrow\{x_4,x_5,x_6\}$ because $d_{[X_6]}^+\geq3$, then we obtain a $4$-ary spanning tree induced by $\{x_{11}x_7,x_{11}x_8,x_{11}x_9,x_{11}u,ux_1,ux_2,ux_3,ux_{10},x_3x_4,x_3x_5,x_3x_6,x_3x_{12}\}$. So we may assume $[X_7]\Rightarrow \{x_{11},x_{12}\}$ and $x_3\Rightarrow [X_2]$ because $d_{[X_7]}^+\geq3$.
Suppose $d_{\{x_4,\dots,x_{10}\}}^+\geq4$, say $x_6\Rightarrow \{x_7,\dots,x_{10}\}$. Then $\{ux_3,\dots,ux_6,x_3x_1,x_3x_2,x_3x_{11},x_3x_{12},x_6x_7,\dots,x_6x_{10}\}$ induces a desired $4$-ary spanning tree. Next we consider the case when $d_{\{x_4,\dots,x_{10}\}}^+=3$.
Since $d^+(x_3)\leq10$, there exists $v\in \{x_4,\dots,x_{10}\}$ such that $v\to x_3$, without loss of generality, say $x_7\Rightarrow \{x_3,\dots,x_6\}$. Then we get a $4$-ary spanning tree induced by $\{ux_7,ux_8,ux_9,ux_{10},x_7x_3,\dots,x_7x_6,x_3x_1,x_3x_2,\\x_3x_{11},x_3x_{12}\}$.
\medskip

{\bf Case 3:} $d^+(u)=9$.

Let $T[\{x_{12},x_{11},x_{10}\}]$ be a transitive $3$-tournament with $x_{12}\Rightarrow \{x_{11},x_{10}\}$ and $x_{11}\to x_{10}$. By Claim \ref{234}, we may assume $[X_9]\Rightarrow x_{12}$. If $d_{[X_9]}^+(x_{11})\geq2$, say $x_{11}\Rightarrow\{x_8,x_9\}$ and $x_4\Rightarrow\{x_5,x_6,x_7\}$ since $d_{[X_7]}^+(x_4)\geq3$. Then we obtain a $4$-ary spanning tree induced by $\{x_{11}x_8,x_{11}x_9,x_{11}x_{10},x_{11}u,ux_1,\dots,ux_4,x_4x_5,x_4x_6,x_4x_7,x_4x_{12}\}$. So we assume $[X_8]\Rightarrow \{x_{12},x_{11}\}$. If $d_{[X_8]}^+(x_{10})\geq3$, say $x_{10}\Rightarrow\{x_6,x_7,x_8\}$ and $x_3\Rightarrow\{x_1,x_2\}$ because $d_{[X_5]}^+\geq2$. Then $\{x_{10}x_6,x_{10}x_7,x_{10}x_8,x_{10}u,ux_3,ux_4,ux_5,ux_9,x_3x_1,x_3x_2,x_3x_{11},x_3x_{12}\}$ induces a desired $4$-ary spanning tree. So we may further assume $[X_6]\Rightarrow \{x_{12},x_{11},x_{10}\}$ and $x_1\to x_2$. Finally, we obtain a $4$-ary spanning tree of $T$ by a similar discussion for $d_{\{x_3,\dots,x_9\}}^+$ as Case 2.

Let $x_{12}\to x_{11}\to x_{10}\to x_{12}$. Suppose $d_{[X_9]}^+(x_{12})\geq3$, say $x_{12}\Rightarrow\{x_7,x_8,x_9\}$. If $d_{[X_6]}^+(x_{11})\geq2$ or $d_{[X_6]}^+(x_{10})\geq2$, say $x_{11}\Rightarrow\{x_5,x_6\}$ or $x_{10}\Rightarrow\{x_5,x_6\}$, then we obtain a $4$-ary spanning tree induced by $\{x_{12}x_7,x_{12}x_8,x_{12}x_9,x_{12}x_{11},x_{11}x_5,x_{11}x_6,x_{11}x_{10},x_{11}u,ux_1,\dots,ux_4\}$ or $\{x_{10}x_5,x_{10}x_6,x_{10}x_{12},x_{10}u,ux_1,\dots,ux_4,x_{12}x_7,x_{12}x_8,x_{12}x_9,x_{12}x_{11}\}$.
So we assume $[X_2]\\\Rightarrow\{x_{12},x_{11},x_{10}\}$ when $d_{[X_9]}^+(x_{12})\leq5$. When $d_{[X_9]}^+(x_{12})\geq6$, we assume $x_{12}\Rightarrow\{x_4,\dots,x_9\}$ and $x_7\Rightarrow\{x_8,x_9\}$ because $d^+_{\{x_4,\dots,x_9\}}\geq2$.
If $x_7\Rightarrow\{x_{11},x_{10}\}$, then we obtain a $4$-ary spanning tree induced by $\{x_{12}x_5,x_{12}x_6,x_{12}x_7,x_{12}u,x_7x_8,\dots,x_7x_{11},ux_1,\dots,ux_4\}$. Since there are at least two vertices with out-degree more than one in $\{x_4,\dots,x_9\}$, say $x_6$ and $x_7$. So we assume $x_{11}\to x_7$, $x_{10}\to x_6$ and $[X_2]\Rightarrow\{x_{12},x_{11},x_{10}\}$. By the symmetry of $x_{12},x_{11}$ and $x_{10}$,
 we get $[X_2]\Rightarrow\{x_{12},x_{11},x_{10}\}$ and $x_1\to x_2$ in each case. Finally, we obtain a $4$-ary spanning tree of $T$ by a similar discussion for $d_{\{x_3,\dots,x_9\}}^+$ as Case 2.\medskip

{\bf Case 4:} $d^+(u)=8$.

By Claim \ref{234}, we may assume $d^+_{N^-(u)}\leq2$. Let $x_{12}\Rightarrow\{x_{11},x_{10}\}$, $x_{11}\Rightarrow\{x_{10},x_9\}$ and $x_9\to x_{12}$.

First we suppose $d_{[X_8]}^+(x_{12})\geq2$, say $x_{12}\Rightarrow\{x_7,x_8\}$. If $d_{[X_6]}^+(x_{11})\geq2$ or $d_{[X_6]}^+(x_9)\geq2$, say
$x_{11}\Rightarrow \{x_5,x_6\}$ or $x_9\Rightarrow \{x_5,x_6\}$, then we get a desired set $\{x_{12}x_{11},x_{12}x_{10},x_{12}x_8,x_{12}x_7,x_{11}x_9,\\x_{11}x_5,x_{11}x_6,x_{11}u,ux_1,\dots,ux_4\}$ or $\{x_9x_{12},x_9x_5,x_9x_6,x_9u,x_{12}x_{11},x_{12}x_{10},x_{12}x_8,x_{12}x_7,ux_1,\dots,\\ux_4\}$. In particular, if $d_{[X_8]}^+(x_{12})\geq3$, say $x_{12}\Rightarrow\{x_6,x_7,x_8\}$, and $d_{[X_5]}^+(x_{11})\geq1$, then we get a $4$-ary spanning tree induced by $\{x_{12}x_{11},x_{12}x_8,x_{12}x_7,x_{12}x_6,x_{11}x_{10},x_{11}x_9,x_{11}x_5,x_{11}u,ux_1,\dots,\\ux_4\}$. Since $d_{[X_8]}^+(x_{12})\leq 5$, we may assume $[X_2]\Rightarrow \{x_{12}, x_{11}, x_9\}$ and $x_1\to x_2$. And it follows that, when $d_{[X_8]}^+(x_{12})\geq2$, $\{x_1x_2, x_1x_{12}, x_1x_{11},x_1x_9, x_{12}x_{10}, x_{12}x_8,x_{12}x_7,x_{12}u,ux_3,\dots,ux_6\}$ induces a desired spanning tree.

 We next consider the case when $d_{[X_8]}^+(x_{12})\leq1$, say $N^+_{[X_8]}(x_{12})\subseteq \{x_8\}$. If $x_9\to x_{10}$, then we assume $[X_5]\Rightarrow \{x_{12}, x_{11}, x_9\}$ by the symmetry of $x_{12},x_{11}$ and $x_9$. If $x_{10}\Rightarrow[X_5]$, say $x_2\to x_1$, then $\{x_{10}x_2,x_{10}x_3,x_{10}x_4,x_{10}u, ux_5,\dots,ux_8, x_2x_1,x_2x_{12},x_2x_{11},x_2x_9\}$ induces a desired $4$-ary spanning tree. So we may further suppose $x_{10}\to x_9$.
If $d_{[X_7]}^+(x_{11})\geq1$, say $x_{11}\to x_7$ and assume $x_4\Rightarrow[X_3]$ because $d^+_{[X_6]}\geq3$, then we obtain a $4$-ary spanning tree induced by $\{x_{11}x_{10},x_{11}x_9,x_{11}x_7,x_{11}u,ux_4,ux_5,ux_6,ux_8,x_4x_1,x_4x_2,x_4x_3,x_4x_{12}\}$. So we may assume $[X_7]\Rightarrow x_{11}$. If $d_{[X_7]}^+(x_{10})\geq2$, say $x_{10}\Rightarrow\{x_6, x_7\}$, assume  $x_3\Rightarrow[X_2]$ because $d_{[X_5]}^+\geq2$, then we obtain a desired set $\{x_{10}x_9,x_{10}x_6,x_{10}x_7,x_{10}u,ux_3,ux_4,ux_5,ux_8,x_3x_1,x_3x_2,x_3x_{11},x_3x_{12}\}$. So we may assume $[X_6]\Rightarrow\{x_{12},x_{11},x_{10}\}$. If $x_9\Rightarrow[X_6]$, say $x_2\to x_1$, then we obtain a $4$-ary spanning tree induced by $\{x_9x_2,x_9x_3,x_9x_4,x_9u,ux_5,\dots,ux_8,x_2x_1,x_2x_{12},x_2x_{11},x_2x_{10}\}$. Consequently, there exists a vertex $v\in[X_6]$ such that $v\Rightarrow\{x_9,\dots,x_{12}\}$, say $v=x_1$. Then we obtain a $4$-ary spanning tree of $T$ by a similar discussion for $d_{\{x_2,\dots,x_8\}}^+$ as Case 2.\medskip

{\bf Case 5:} $d^+(u)=7$.

Suppose $d^+_{N^-(u)}=d^+_{N^-(u)}(x_{12})$.

Firstly, suppose $d^+_{N^-(u)}(x_{12})=4$, say $x_{12}\Rightarrow\{x_8,\dots,x_{11}\}$.  If there exists some vertex, say $x_4$, such that $d^+_{[X_7]}(x_4)\geq3$ and $x_4\to x_{12}$, then we assume $x_4\Rightarrow\{x_5,x_6,x_7\}$ and obtain a $4$-ary spanning tree induced by $\{ux_1,\dots,ux_4,x_4x_5,x_4x_6,x_4x_7,x_4x_{12},x_{12}x_8,\dots,x_{12}x_{11}\}$. Since $d^+(x_{12})\leq7$, we may assume $x_{12}\Rightarrow\{x_6,x_7\}$ and $T[[X_5]]$ is $2$-regular with $x_1\Rightarrow \{x_2,x_3\}$. Let $d^+_{N^-(u)}(x_{11})\geq2$. Since $d^+(x_{11})\leq7$, there exists some vertex $v\in[X_5]$ such that $v\to x_{11}$, say $v=x_1$. Then  $\{x_1x_2,x_1x_3,x_1x_{11},x_1x_{12},x_{12}x_8,x_{12}x_9,x_{12}x_{10},x_{12}u,ux_4,\dots,ux_7\}$ induces a desired $4$-ary spanning tree.

Next, suppose $d^+_{N^-(u)}(x_{12})=3$, say $x_{12}\Rightarrow\{x_{11},x_{10},x_9\}$. If there exists some vertex, say $x_4$, such that $d^+_{[X_7]}(x_4)\geq3$ and $x_4\to x_8$, then we assume $x_4\Rightarrow\{x_5,x_6,x_7\}$ and obtain a $4$-ary spanning tree induced by $\{x_{12}x_{11},x_{12}x_{10},x_{12}x_9,x_{12}u,ux_1,\dots,ux_4,x_4x_5,x_4x_6,x_4x_7,x_4x_8\}$.
If $d^+_{[X_7]}(x_8)\geq3$, say $x_8\Rightarrow\{x_5,x_6,x_7\}$, then we get a $4$-ary spanning tree induced by $\{x_8x_5,x_8x_6,\\x_8x_7,x_8x_{12},x_{12}x_{11},x_{12}x_{10},x_{12}x_9,x_{12}u,ux_1,\dots,ux_4\}$. So we may assume $x_8\Rightarrow\{x_6,x_7\}$, $T[[X_5]]$ is $2$-regular with $x_1\Rightarrow \{x_2,x_3\}$ and $\{x_6,x_7\}\Rightarrow[X_5]$. If $d^+_{[X_5]}(x_{12})\geq1$, say $x_{12}\to x_5$, then we obtain a $4$-ary spanning tree induced by $\{x_8u,x_8x_6,x_8x_7,x_8x_{12},,x_6x_1,\dots,x_6x_4,x_{12}x_{11},\\x_{12}x_{10},x_{12}x_9,x_{12}x_5\}$. So we may assume $[X_5]\Rightarrow \{x_8,x_{12}\}$. Then  we obtain a $4$-ary spanning tree induced by $\{x_1x_2,x_1x_3,x_1x_8,x_1x_{12},x_{12}x_9,x_{12}x_{10},x_{12}x_{11},x_{12}u,ux_4,\dots,ux_7\}$.

Finally, we consider the case when $T[{N^-(u)}]$ is $2$-regular, say $x_{12}\Rightarrow\{x_{11},x_{10}\}$ and $x_{11}\to x_{10}$. If $d^+(v)\leq6$ for any $v\in V(T)\less\{u\}$, then the out-degree sequence of $T$ is $\{5,6,\dots,6,7\}$. So there exist three vertices, say $x_5$, $x_6$ and $x_7$, such that $x_{10}\Rightarrow\{x_6,x_7\}$ and $x_{12}\to x_5$. Then we obtain a desired  tree induced by $\{x_{12}x_{11},x_{12}x_{10},x_{12}x_5,x_{12}u,x_{10}x_6,\dots,x_{10}x_9,ux_1,\dots,ux_4\}$. We next consider the remaining two cases. If $d^+(x_{12})=7$, say
$x_{12}\Rightarrow[X_4]$, we may assume $T[\{x_5,\dots,x_9\}]$ is $2$-regular with $x_9\Rightarrow \{x_8,x_7\}$ by the symmetry of $x_{12}$ and $u$, then we obtain a desired tree induced by $\{x_9x_8,x_9x_7,x_9x_{12},x_9u,x_{12}x_{11},x_{12}x_{10},x_{12}x_1,x_{12}x_2,ux_3,\dots,ux_6\}$. Without loss of generality, if $d^+(x_1)=7$, we may assume $x_1\Rightarrow\{x_4,\dots,x_{10}\}$ because $T[N^-(x_1)]$ is $2$-regular by the symmetry of $x_1$ and $u$, then we obtain a $4$-ary spanning tree induced by $\{x_{12}x_{11},x_{12}x_{10},x_{12}x_1,x_{12}u,x_1x_6,\dots,x_1x_9,ux_2,\dots,ux_5\}$.\medskip

{\bf Case 6:} $d^+(u)=6$.

 In this case, $T$ is $6$-regular. Firstly, suppose $d^+_{N^-(u)}=d^+_{N^-(u)}(x_{12})=5$, say $x_{12}\Rightarrow\{x_7,\dots,x_{11}\}$. Let $1\leq d^+_{N^-(u)}(x_8)\leq3$ and assume $x_8\Rightarrow\{x_5,x_6,x_7\}$. Then we obtain a $4$-ary spanning tree induced by $\{x_{12}x_{11},\dots,x_{12}x_8,x_8x_7,x_8x_6,x_8x_5,x_8u,ux_1,\dots,ux_4\}$.
Then, suppose $d^+_{N^-(u)}=4$, say $x_{12}\Rightarrow\{x_8,\dots,x_{11}\}$. If $d_{[X_6]}^+(x_7)\geq2$, say $x_7\Rightarrow\{x_5,x_6\}$, then we obtain a desired set $\{x_7x_5,x_7x_6,x_7x_{12},x_7u,x_{12}x_{11},\dots,x_{12}x_8,ux_1,\dots,ux_4\}$. Notice that $T$ is $6$-regular, so we may assume $x_7\Rightarrow\{x_6,x_8,x_9,x_{10}\}$. If $d_{[X_5]}^+(x_{12})\geq1$, say $x_{12}\to x_5$, then we obtain a desired tree induced by $\{x_7x_6,x_7x_8,x_7x_{12},x_7u,x_{12}x_{11},x_{12}x_{10},x_{12}x_9,x_{12}x_5,ux_1,\dots,ux_4\}$. So we may assume $[X_5]\Rightarrow x_{12}$ and $x_3\Rightarrow\{x_1,x_2\}$ because $d_{[X_5]}^+\geq2$, and then we obtain a $4$-ary spanning tree induced by $\{ux_3,\dots,ux_6,x_3x_1,x_3x_2,x_3x_7,x_3x_{12},x_{12}x_{11},\dots,x_{12}x_8\}$.
Finally, suppose $d^+_{N^-(u)}=3$, say $x_{12}\Rightarrow\{x_9,x_{10},x_{11}\}$ and $x_7\to x_8$. Since $d_{[X_6]}^+(x_7)\geq2$, we assume  $x_7\Rightarrow\{x_5,x_6\}$. Then  $\{x_7x_5,x_7x_6,x_7x_8,x_7x_{12},x_{12}x_{11},x_{12}x_{10},x_{12}x_9,x_{12}u,ux_1,\dots,ux_4\}$ induces a desired $4$-ary spanning tree.

Now suppose $k=5$.
Let $T_{12}$ be a $12$-tournament with $V(T_{12})=\{0,1,\dots,11\}$ and $E(T_{12})=\{(0,3),(0,5),(0,9),(0,10),(0,11),(1,0),(1,4),(1,6),(1,8),(1,9),(1,11),(2,0),(2,1),\\(2,7),(2,8),
(2,10),(2,11),(3,1),(3,2),(3,6),(3,9),(3,10),(4,0),(4,2),(4,3),(4,7),(4,9),\\(5,1),(5,2),(5,3),(5,4),(5,8),
(5,11),(6,0),(6,2),(6,4),(6,5),(6,10),(7,0),(7,1),(7,3),(7,5),\\(7,6),(8,0),(8,3),(8,4),(8,6),(8,7),(9,2),(9,5),(9,6),(9,7),(9,8),
(9,11),(10,1),(10,4),\\(10,5),(10,7),(10,8),(10,9),(11,3),(11,4),(11,6),
(11,7),(11,8),(11,10)\}$. It is easy to check that $T_{12}$ satisfies the condition of Observation \ref{obs}. Therefore, $T_{12}$ contains no $5$-ary spanning tree, which implies that $h(5)\geq13$.\qed

\begin{rem}
Using the similar method as $h(4)$, we can prove that $h(5)=13$. However, the proof is too long to include here. Some new methods are needed to determine the exact values of $h(k)$ for $k\geq 5$.
\end{rem}

\section{Acknowledgments} The authors would like to thank the two referees and the editor for their valuable suggestions and comments. This work was partially supported by the National Natural Science Foundation of China (No. 11922112) and the Fundamental Research Funds for the Central Universities.

\end{document}